\documentclass[conference]{IEEEtran}
\IEEEoverridecommandlockouts
\usepackage{cite}

\usepackage[utf8]{inputenc}
\usepackage{amsmath,amssymb,amsfonts}
\usepackage{algorithmic}
\usepackage{graphicx}
\usepackage{textcomp}
\usepackage{xcolor}
\usepackage{commath} 

\def\BibTeX{{\rm B\kern-.05em{\sc i\kern-.025em b}\kern-.08em
    T\kern-.1667em\lower.7ex\hbox{E}\kern-.125emX}}

\newtheorem{theorem}{Theorem}

\newtheorem{lemma}{Lemma}

\begin{document}

\title{The dual approach to non-negative super-resolution: impact on primal reconstruction accuracy
\thanks{This publication is based on work supported by the 
EPSRC Centre For Doctoral Training in Industrially Focused 
Mathematical Modelling (EP/L015803/1) in collaboration with 
the National Physical Laboratory and by the Alan Turing Institute 
under the EPSRC grant EP/N510129/1 and the Turing Seed Funding grant SF019.}
}

\author{\IEEEauthorblockN{St\'{e}phane Chr\'{e}tien}
\IEEEauthorblockA{
\textit{National Physical Laboratory}\\
London, UK \\
stephane.chretien@npl.co.uk}
\and
\IEEEauthorblockN{Andrew Thompson}
\IEEEauthorblockA{
\textit{National Physical Laboratory}\\
London, UK \\
andrew.thompson@npl.co.uk}
\and
\IEEEauthorblockN{Bogdan Toader}
\IEEEauthorblockA{
\textit{University of Oxford}\\
Oxford, UK \\
toader@maths.ox.ac.uk}
}

\maketitle

\begin{abstract}
We study the problem of super-resolution, where we 
recover the locations and weights of non-negative point sources
from a few samples of their convolution with a Gaussian kernel. 
It has been recently shown that exact recovery is possible 
by minimising the total variation norm of the measure. 
An alternative practical approach is to solve its dual. In this paper, 
we study the stability of solutions with respect to the solutions 
to the dual problem.
In particular, we establish a relationship between perturbations in the dual variable and the primal variables around the
optimiser. This is achieved by applying a quantitative
version of the implicit function theorem in a non-trivial way.
\end{abstract}


\section{Problem setup}

In the study of non-negative super-resolution, we want
to estimate a signal $x$ which consists of a number of
point sources with unknown locations and non-negative magnitudes,
from only a few measurements of the convolution of $x$ with
a known kernel $\phi$. This is a problem that 
arises in a number of applications, for example 
fluorescence microscopy~\cite{betzig2006imaging},
astronomy~\cite{puschmann2005super}
and ultrasound imaging~\cite{tur2011innovation}.
In such applications, the measurement
device has limited resolution and cannot distinguish 
between distinct point sources in the input signal $x$. 
This is often modelled as a deconvolution problem with a Gaussian kernel.

Specifically, let $x$ be a 
non-negative measure on $I=[0,1]$ consisting
of $K$ unknown non-negative 
point sources:
\begin{equation*}
  x = \sum_{k=1}^{K} a_k \delta_{t_k},
\end{equation*}
with $a_k > 0$, for all $k=1,\ldots,K$, 
and let $\mathbf{y} \in \mathbb{R}^M$ be
the vector of measurements
obtained by sampling 
the convolution of $x$ with a known 
kernel $\phi$ (e.g. a Gaussian) 
at locations $s_j$:
\begin{equation*}
    y_j =  
    \int_I \phi(t-s_j) x(\dif t),
\end{equation*} 
for all $j=1,\ldots,M$.
Note that, because $x$ is a discrete measure,
each entry in $\mathbf{y}$ is of the form:
\begin{equation*}
   y_j =  \sum_{k=1}^{K} a_k \phi(t_k-s_j),
\end{equation*}
for all $j = 1,\ldots, M$. Let 
$
    \Phi(t) = [\phi(t-s_1),\ldots,\phi(t-s_M)]^T,
$
then $x$ can be recovered by solving the 
following program:
\footnote{
    We assume that the measurements
    are exact. We treat the case when
    the measurements are noisy in the journal version of 
    this paper.
}
\begin{equation}
  \min_{\hat{x} \geq 0} \int_I \hat{x}(\dif t)
  \quad \text{subject to} \quad
  \mathbf{y} = \int_I \Phi(t) \hat{x}(\dif t).
  \label{eq:main prog}
\end{equation}
The problem of super-resolution has been studied extensively
in the literature since the seminal paper \cite{candes2014towards}, 
which addressed the case of complex amplitudes. 
In \cite{schiebinger}, the authors showed that when the coefficients $a_k$,
$k=1,\ldots,K$ are positive, which corresponds to the present setting, 
exact recovery is possible without separation of sources.

The dual of problem \eqref{eq:main prog} is 
\begin{equation}
  \max_{\lambda \in \mathbb{R}^M}\;
  \mathbf{y}^T \mathbf{\lambda}
  \quad \text{subject to} \quad 
  \sum_{j=1}^M \lambda_j \phi(t - s_j) \leq 1, 
  \quad \forall t \in I.
  \label{eq:dual}
\end{equation}
The dual problem \eqref{eq:dual} is a finite-dimensional problem 
with infinitely many constraints, known as a semi-infinite program. 
Such problems can be solved using a number of algorithms including 
exchange methods~\cite{eftekhari2019equivalence} and sequential 
quadratic programming~\cite{lopez2007semiinfinite}. The advantage
over algorithms that solve the primal problem 
(for example the ADCG algorithm~\cite{boyd2017}) 
is working in a finite dimensional space, which
simplifies the analysis.

Consider a solution $\lambda^*$ of the dual 
problem \eqref{eq:dual} which corresponds to a dual 
certificate, namely a function
\begin{equation}
  q(s) = \sum_{j=1}^{M} \lambda_j^* \phi(s-s_j),
  \label{eq:dual cert}
\end{equation}
which satisfies the conditions:
\begin{align}
    &q(t_i) = 1, \quad \forall i=1,\ldots,K,
    \label{eq:cond1}
    \\
    &q(s) < 1, \quad \forall s \ne t_i, 
    \quad \forall i=1,\ldots,K.
    \label{eq:cond2}
\end{align}
Then the local maximisers of $q(s)$ correspond to the 
source locations $\{t_k\}_{k=1}^K$ and the 
amplitudes $\{a_k\}_{k=1}^K$ are found
by solving a linear system.

In this paper, we analyse how small perturbations 
of $\lambda^*$ affect the local maximisers
of $q(s)$ in the case when the convolution kernel
is Gaussian $\phi(t) = e^{-t^2/\sigma^2}$.
The outcome is a bound on how far the 
estimated locations $t_k$ and magnitudes $a_k$ are from
their true values obtained for exact $\lambda^*$.
This gives us an insight
into the size of the error in the locations and magnitudes
when we apply an optimisation algorithm to the
dual of the super-resolution problem.

\section{Bound on the error as $\lambda$ is perturbed}

In this section we present our main results, namely
two theorems that give bounds on the perturbations
around the source locations $t_k$ and the 
magnitudes $a_k$ respectively, as the dual variable
is perturbed away from the optimiser $\lambda^*$, 
when the convolution kernel is a 
Gaussian with known width $\sigma$:
$
    \phi(t) = e^{-t^2/\sigma^2}.
$

\begin{theorem}
  \textbf{(Dependence of $t$ on $\lambda$)}
  Let $\lambda^* \in \mathbb{R}^M$ be a solution of the 
  dual program \eqref{eq:dual} with $\phi$ Gaussian 
  as given above such that the dual
  certificate $q(s)$ defined in \eqref{eq:dual cert}
  satisfies conditions \eqref{eq:cond1} and \eqref{eq:cond2},
  $\lambda$ a perturbation of $\lambda^*$ in a ball of 
  radius $\delta_{\lambda}$ and $t$ an arbitrary local
  maximiser of $q_{\lambda}(s)=\sum_{j=1}^M \lambda_j \phi(s-s_j)$ 
  so that for $\lambda=\lambda^*$, 
  the corresponding local maximiser $t^*$ is a source 
  location of $x$. Then
  \begin{equation*}
    |t-t^*| \leq C_{t^*} \|\lambda-\lambda^*\|_2,
  \end{equation*}
  provided that
  \begin{equation*}
    \delta_{\lambda} \leq
    \frac{
      |q''(t^*)|^2 \sigma^3 \sqrt{e}
    }{
      4\sqrt{2} \left(
        2 + c 
        R 
      \right)
      M
    },
    \label{eq:delta_lambda}
  \end{equation*}
  where
  \begin{align}
    C_{t^*} &= \frac{1}{4 + cR}
      \left[
        1 + \frac{
          2\sqrt{2M} (2+cR)
        }{
          |q''(t^*)| \sqrt{e} 
        }
      \right],
    \label{eq:C_t}
    \\
    R &= \frac{ \| \lambda^* \|_2}{\sigma},
  \end{align}
  and $c \approx 3.9036$ is a universal constant.
  \label{thm:t dep lambda}
\end{theorem}

One of the main conclusions which can be drawn from this result 
is that the primal spike location error is controlled 
in $l_\infty$, but degrades as a function of the number of measurements 
in the order of $\sqrt{M}$. Of crucial importance is the curvature 
of the dual certificate at the true solution: the flatter the 
certificate, the worse the estimation error. 
Our theorem also gives important information about the accuracy 
in the dual variable required to guarantee our upper bound on the 
error of recovery. This accuracy is of the inverse order of the 
number of measurements, which is quite a stringent constraint. 
Both the $M$ and the $\sqrt{M}$ factors are a consequence of the
way we bound sums of shifted copies of the kernel, namely
$\sum_{j=1}^M \phi(t-s_j) \leq M \max_{t \in \mathbb{R}} \phi(t)$.
Given the fast decay of the Gaussian, it is clear that this is not
a tight bound. However, any bound would reflect the density
of samples close to each source location.

We will now give a result regarding the perturbation
of the magnitudes $a_k$ when $\lambda^*$ is perturbed.
Let $\Phi$ be the matrix whose entries are defined as
\begin{equation}
    \Phi_{ij} = \phi(t_j-s_i),
    \label{eq:def Phi}
\end{equation}
and $\mathbf{t^*}$ and $\mathbf{a^*}$ the vectors of 
true source locations and weights:
\begin{align*}   
    \mathbf{t^*} = [t_1,\ldots,t_K]^T,
    \quad\quad
    \mathbf{a^*} = [a_1,\ldots,a_K]^T.
\end{align*}
When we solve \eqref{eq:dual} exactly, we obtain
the source locations by finding the local maximisers
of $q(s)$. Then the vector of weights $\mathbf{a^*}$ 
is found by solving the system
\begin{equation*}
    \Phi \mathbf{a} = \mathbf{y}.
\end{equation*}
When the source locations are perturbed, 
we denote the resulting perturbed data matrix by:
\begin{equation}
  \tilde{\Phi} = \Phi + E,
  \label{eq:phi tilde def}
\end{equation}
and we calculate the vector of perturbed weights
$\tilde{\mathbf{a}}$ as the solution of the least 
squares problem
\begin{equation}
  \min_{\mathbf{a}} \| \tilde{\Phi} \mathbf{a} - \mathbf{y} \|_2.
  \label{eq:perturbed phi least sq}
\end{equation}
The following theorem gives us a bound on the difference
$\| \mathbf{a^*} - \tilde{\mathbf{a}} \|_2$ 
between the vector of true weights $\mathbf{a^*}$
and the vector of weights $\tilde{\mathbf{a}}$ obtained by solving
the least squares problem  \eqref{eq:perturbed phi least sq}
with the perturbed matrix $\tilde{\Phi}$, as a function
of the difference $\|\tilde{\mathbf{t}}-\mathbf{t^*}\|_2$ 
between the perturbed source 
locations $\tilde{\mathbf{t}}$ and the true source 
locations $\mathbf{t^*}$.

\begin{theorem}
  \textbf{(Dependence of $a$ on $t$)}
  Let $\mathbf{t^*} \in [0,1]^{K}$ be the vector of true source 
  locations and $\tilde{\mathbf{t}} \in [0,1]^K$ the perturbed source 
  locations, such that:
    \begin{equation}
        \|\mathbf{t^*}-\mathbf{\tilde{t}}\|_2 
        <
        \frac{
            \sigma^2 \sigma_{\max}(\Phi)
        }{
            4e^{4/\sigma^2} \sqrt{M}
        }
        \left(
            \sqrt{
                1 + \frac{
                    \sigma_{\min}^2(\Phi)
                }{
                    \sigma_{\max}^2(\Phi)
                } 
            } -1
        \right).
        \label{eq:cond t}
    \end{equation} 
  Then the error between the true weights $\mathbf{a^*}$
  and the perturbed weights $\tilde{\mathbf{a}}$ obtained by
  solving program \eqref{eq:perturbed phi least sq}
  is bounded by:
  \begin{equation*}
    \|\mathbf{a^*}-\tilde{\mathbf{a}}\|_2 
    \leq
    \frac{
      4 e^{\frac{4}{\sigma^2}}\sqrt{M}
      \|\mathbf{a^*}\|_2
    }{
      \sigma^2 \sigma_{\min}(\Phi)
    } \|\tilde{\mathbf{t}}-\mathbf{t^*}\|_2
    + O(\|\tilde{\mathbf{t}}-\mathbf{t^*}\|_2^2).
  \end{equation*}
  \label{thm:a dep t}
\end{theorem}

\section{Proofs}

In this section we present the proofs of the two theorems.
Due to space limitations, we skip some details, 
which will be present in the journal version of the paper.

\subsection{Proof of Theorem \ref{thm:t dep lambda}}
  \label{sec:proof thm t dep lambda}

Let $t^*$ be an arbitrary local maximiser of the 
function $q(t)$ in \eqref{eq:dual cert}, so $t^*$ 
is also a source location, and $\lambda^*$ 
the solution to \eqref{eq:dual}.
The key step in this proof is applying a 
quantitative version of the 
Implicit Function Theorem \cite{liverani}
to the function:
\begin{equation}\label{eq:Fdef}
  F(t,\lambda) = \sum_{j=1}^{M} \lambda_j \phi'(t-s_j),
\end{equation}
where $F(t^*,\lambda^*) = 0$ because $t^*$ is a maximizer of $q(s)$ in \eqref{eq:dual cert}. The theorem tells us that
we can express $t$ as a function $t(\lambda)$ of $\lambda$ with:
\begin{equation}
    \partial_{\lambda} t(\lambda) = 
    -\left[
        \partial_t F(t(\lambda),\lambda)
    \right]^{-1}
    \partial_{\lambda} F(t(\lambda),\lambda),
    \label{eq:qift res}
\end{equation}
for $t$ in a ball of radius $\delta_0$ around $t^*$ and 
for $\lambda$ in a ball of radius $\delta_1 \leq \delta_0$ 
around $\lambda^*$, where $\delta_0$ is chosen such that
\begin{equation}
  \sup_{(t,\lambda) \in V_{\delta}}
  \left\| 
    I - 
    \left[
      \partial_t F(t^*,\lambda^*)
    \right]^{-1}
    \partial_t F(t,\lambda)
  \right\|
  \leq \frac12,
  \label{eq:cond delta0}
\end{equation}
where
$
  V_{\delta}= \left\{ 
    (t, \lambda) \in \mathbb{R}^{M+1}:
    |t - t^*| \leq \delta_0,
    \|\lambda - \lambda^*\| \leq \delta_0
  \right\}
$  
and $\delta_1$ is given by 
\begin{equation}
    \delta_1 = (2M_t B_{\lambda})^{-1} \delta_0,
    \label{eq:cond delta1}
\end{equation}
where
\begin{align*}
  B_{\lambda} &=
  \sup_{(t,\lambda) \in V_{\lambda}}
    \| \partial_{\lambda} F(t,\lambda) \|_2,
  \\
  M_t &= \left\| 
    \partial_t F(t^*, \lambda^*)^{-1} 
  \right\|_2.
\end{align*}
The following two lemmas give values of $\delta_0$
and $\delta_1$ that define balls around $t^*$ and $\lambda^*$
respectively which are included in the balls required
by the Quantitative Implicit Function Theorem with radii 
defined in \eqref{eq:cond delta0} and \eqref{eq:cond delta1}.
\begin{lemma}
    \textbf{(Radius of ball around $t^*$)}
    The condition \eqref{eq:cond delta0} is satisfied if
    \begin{equation*}
        \delta_0 
        =\frac{
            \sigma^2 |q''(t^*)| 
          }{
            \sqrt{M}\left(
              4 + 2c \cdot
              \frac{\|\lambda^*\|_2}{\sigma}
            \right)
          }.
        \label{eq:apriori delta bound}
    \end{equation*}
    \label{lem:bnd t}
\end{lemma}

\begin{lemma}
    \textbf{(Radius of ball around $\lambda^*$)}
    For $\delta_0$ from Lemma \ref{lem:bnd t} 
    and $\delta_1$ from condition \eqref{eq:cond delta1}, 
    we have that $\delta_{\lambda} < \delta_1$ if
    \begin{equation*}
      \delta_{\lambda} 
      = \frac{
        \sigma \sqrt{e} |q''(t^*)| 
      }{
        2\sqrt{2M}
      } \cdot \delta_0.
    \end{equation*}
    \label{lem:bnd lambda}
\end{lemma}
Given the definition of the function $F$ in \eqref{eq:Fdef}, we have that
\begin{align*}
  \partial_{t} F(t,\lambda) &= \sum_{j=1}^{m} \lambda_j \phi''(t-s_j), 
  \\
  \partial_{\lambda} F(t,\lambda) &= [\phi'(t-s_1), \quad \ldots \quad, \phi'(t-s_M)]^T. 
\end{align*}
By applying a Taylor expansion to $t(\lambda)$ around
$\lambda^*$ in the region defined by $\delta_0$ and $\delta_{\lambda}$,
we have that
\begin{equation*}
  t(\lambda) = t(\lambda^*) + 
  \left<
    \lambda-\lambda^*,
    \partial_{\lambda} t(\lambda_{\delta})
  \right>,
\end{equation*}
for some $\lambda_{\delta}$ on the line segment 
determined by $\lambda^*$ and $\lambda$,
so
\begin{align}
  \left|
    t(\lambda)-t(\lambda^*)
  \right|  
  \leq
  \left\|
    \lambda-\lambda^*
  \right\|_2
  \cdot
  \left\|
    \partial_{\lambda} t(\lambda_{\delta})
  \right\|_2.
  \label{eq:terms_bnd}
\end{align}
Bounding the second factor in \eqref{eq:terms_bnd} above,
using \eqref{eq:qift res} and applying a number
of manipulations, we
obtain
\begin{equation}
  \delta_t \leq
  \frac{C + \delta_t A}{B}
  \cdot \delta_{\lambda}.
  \label{eq:equiv terms bnd}
\end{equation}
where $\delta_t = |t(\lambda)-t(\lambda^*)|$,
$\|\lambda - \lambda^*\|_2 \leq \delta_{\lambda}$ and
\begin{align*}
  A &=
    \left\|
      \left[ 
        \phi''(t(\lambda)-s_j)
      \right]_{j=1}^M
    \right\|_2,
  \\
  B &=
    \left|
      \sum_{j=1}^M 
      \lambda_j^* \phi''(t(\lambda)-s_j)
    \right|,
  \\
  C &= 
    \left\| \left[
      \phi'(t(\lambda_{\delta})-s_j)
    \right]_{j=1}^M \right\|_2.
\end{align*}
We now need to lower bound $B$ and 
upper bound $C+\delta_t A$.
For $B$, by applying a Taylor expansion around
$ t(\lambda^*) - s_j$, each term in the sum is equal to
\begin{align*}
    \lambda_j^* \phi''(t(\lambda^*)-s_j)
    + \left(t(\lambda) - t(\lambda^*)\right)
    \lambda_j^* \phi'''(\xi_j),
\end{align*}
for $j=1,\ldots,M$ and some $\xi_j$ in the interval 
$$
  \left[
    t(\lambda^*) - s_j - |t(\lambda)-t(\lambda^*)|,
    t(\lambda^*) - s_j + |t(\lambda)-t(\lambda^*)|
  \right].
$$
By combining this with the reverse triangle inequality,
Cauchy-Schwartz inequality and Lemma \ref{lem:bnd t}, 
we obtain
\begin{equation*}
  B \geq |q''(t^*)|
  \left[
    1 - \frac{
      c \| \lambda^* \|_2 
    }{
      4 \sigma + 2c \| \lambda^* \|_2
    }
  \right],
  \label{eq:bound B}
\end{equation*}
where $c\approx 3.9036$. By using the global maximum 
of $\phi'(t)$ and $\phi''(t)$ for $\phi(t)$ Gaussian,
we have that
\begin{align*}
    A &\leq \frac{2\sqrt{M}}{\sigma^2}
    \quad\text{and}\quad 
    C \leq \frac{\sqrt{2M}}{\sigma\sqrt{e}}.
\end{align*}
By finally plugging the above inequalities 
and the result of Lemma \ref{lem:bnd t}
into \eqref{eq:equiv terms bnd}, we obtain
the expression for $C_{t^*}$ in \eqref{eq:C_t}.

\subsection{Proof of Theorem \ref{thm:a dep t}}
  \label{sec:proof thm a dep t}

We apply equation (4.2) from \cite{Stewart90perturbationtheory}
with $e = 0$ (the noise in the observations)
and obtain
\begin{equation}
  \mathbf{\tilde{a}} = \mathbf{a^*} 
    - \Phi^{\dagger} E\mathbf{a^*} - F^T E \mathbf{a^*},
  \label{eq:a tilde}
\end{equation}
where $\Phi^{\dagger} = (\Phi^T \Phi)^{-1} \Phi^T$ is the
pseudoinverse of $\Phi$ and $F=O(E)$ is the perturbation
of the $\Phi^{\dagger}$ due to the perturbation $E$ 
of $\Phi$, namely
\begin{equation*}
  \tilde{\Phi}^{\dagger} = \Phi^{\dagger} + F^T.
\end{equation*}
In order to obtain an explicit expression for $F$, 
we write $\tilde{\Phi}^{\dagger}$:
\begin{align}
  \tilde{\Phi}^{\dagger} 
  &= (\tilde{\Phi}^T \tilde{\Phi})^{-1} \tilde{\Phi}^T 
  \nonumber \\
  &= \left[ (\Phi+E)^T (\Phi+E) \right]^{-1} (\Phi+E)^T
  \quad\quad\text{by } \eqref{eq:phi tilde def}
  \nonumber \\
  &= (\Phi^T \Phi + \Delta)^{-1} (\Phi^T + E^T),
  \label{eq:tilde phi expl}
\end{align}
where 
\begin{equation}
  \Delta = E^T \Phi + \Phi^T E + E^T E
  \in \mathbb{R}^{K \times K}.
  \label{eq:delta def}
\end{equation}
We then write the first factor in \eqref{eq:tilde phi expl} 
as
\begin{align}
  (\Phi^T \Phi + \Delta)^{-1} 
  &= \left[
    \Phi^T \left(
      I + {\Phi^{\dagger}}^T \Delta \Phi^{\dagger} 
    \right) \Phi
  \right]^{-1}
  \nonumber \\
  &= \Phi^{\dagger} \left[
    I + \sum_{k=1}^{\infty} (-1)^k 
    \left( {\Phi^{\dagger}}^T \Delta \Phi^{\dagger} \right)^k
  \right]  {\Phi^{\dagger}}^{T} 
  \nonumber \\
  &= (\Phi^T \Phi)^{-1} + S_{\Phi},
  \label{eq:with Neumann}
\end{align}
where 
\begin{equation}
  S_{\Phi} = 
  \Phi^{\dagger} \left[
    \sum_{k=1}^{\infty} (-1)^k 
    \left( {\Phi^{\dagger}}^T \Delta \Phi^{\dagger} \right)^k
  \right]  {\Phi^{\dagger}}^{T}
  \in \mathbb{R}^{K \times K},
    \label{eq:Amatrix}
\end{equation}
and in the second inequality in \eqref{eq:with Neumann}
we applied the Neumann series expansion to the 
matrix $-{\Phi^{\dagger}}^T \Delta \Phi^{\dagger}$, 
which converges if 
\begin{equation}
    \|
        -{\Phi^{\dagger}}^T \Delta \Phi^{\dagger}       
    \|_2 < 1.
    \label{eq:cond Neumann}
\end{equation}
We will return to condition \eqref{eq:cond Neumann} 
at the end of this section.
We now substitute \eqref{eq:with Neumann} 
in \eqref{eq:tilde phi expl}, giving
\begin{align*}
  \tilde{\Phi}^{\dagger}
  &= \left[
    (\Phi^T \Phi)^{-1} + S_{\Phi}
  \right](\Phi^T + E^T)
  \nonumber \\
  &= \Phi^{\dagger} + (\Phi^T \Phi)^{-1} E^T 
    + S_{\Phi} \Phi^T + S_{\Phi} E^T,
\end{align*}
so we have that
\begin{equation}
  F^T = (\Phi^T \Phi)^{-1} E^T + S_{\Phi} \Phi^T + S_{\Phi} E^T,
  \label{eq:F expl}
\end{equation}
which is indeed $O(E)$, since $S_{\Phi} = O(\Delta)$ 
and $\Delta=O(E)$.
We next upper bound $\|S_{\Phi}\|_2$. 
From \eqref{eq:Amatrix} we have
\begin{equation}
  \| S_{\Phi} \|_2 \leq
  \| \Phi^{\dagger} \|^2_2
  \sum_{k=1}^{\infty}
  \| \Phi^{\dagger}\|^{2k}_2 \| \Delta \|^k_2.
  \label{eq:A first upper bnd}
\end{equation}
Now let $D$ be an upper bound on $\|\Delta\|_2$, obtained by
applying the triangle inequality in \eqref{eq:delta def}, so that
\begin{equation}
  \| \Delta \|_2 \leq D
  = 2 \| E \|_2 \| \Phi\|_2 + \|E\|^2_2.
  \label{eq:def D}
\end{equation}
Then, from \eqref{eq:A first upper bnd} we have
\begin{align}
  \|S_{\Phi}\|_2 &\leq 
  \| \Phi^{\dagger} \|^2_2
  \sum_{k=1}^{\infty}
  \| \Phi^{\dagger}\|^{2k}_2 D^k
  = \frac{
    D \| \Phi^{\dagger}\|^4_2
  }{
    1 - D \| \Phi^{\dagger}\|^2_2
  },
  \label{eq:bound norm A}
\end{align}  
where the series converges if $D\|\Phi^{\dagger}\|_2^2 < 1$,
in which case the denominator in the last fraction above is 
positive. We return to this condition at the end of 
the section.
We also know that\footnote{
  Using the SVD $\Phi = U \Sigma V^T$, we have
  $  
    \Phi^{\dagger} = (\Phi^T \Phi)^{-1} \Phi^T
    = (V \Sigma^2 V^T)^{-1} V \Sigma U^T 
    = V \Sigma^{-1} U^T
  $,
  so the conclusion follows.
}
\begin{equation}
  \| \Phi^{\dagger} \|_2 = \frac{1}{\sigma_{\min}(\Phi)}.
  \label{eq:norm phi dagger}
\end{equation}
By applying the triangle inequality in \eqref{eq:F expl} and 
then using \eqref{eq:bound norm A} and the fact 
that $\|(\Phi^T\Phi)^{-1}\|_2 = 1/\sigma_{\min}^2(\Phi)
=\|\Phi^{\dagger}\|_2^2$ (from \eqref{eq:norm phi dagger}),
we obtain
\begin{equation}
  \|F\|_2 \leq 
  \|E\|_2  \|\Phi^{\dagger}\|_2^2
  + \frac{
    D \| \Phi^{\dagger}\|^4_2
  }{
    1 - D \| \Phi^{\dagger}\|^2_2
  } 
  \left( \| \Phi \|_2 + \|E\|_2 \right),
  \label{eq:norm F}
\end{equation}
where $D$ is given in \eqref{eq:def D}. 
It remains to establish an upper bound 
on $\|E\|_F$, and consequently on $\|E\|_2$.
The following lemma gives us such a bound.

\begin{lemma}
    \textbf{(Upper bound of $\|E\|_F)$}
    Let $E = \tilde{\Phi} - \Phi$ for $\Phi$ and $\tilde{\Phi}$ 
    as defined in \eqref{eq:def Phi} and \eqref{eq:phi tilde def} 
    respectively for $t_j, \tilde{t}_j \in [0,1]$ 
    for $j=1,\ldots,K$. Then:
    \begin{equation}
        \|E\|_F 
        \leq
        \frac{4 e^{\frac{4}{\sigma^2}}\sqrt{M}}{\sigma^2} 
        \|\mathbf{\tilde{t}}-\mathbf{t^*}\|_2.
        \label{eq:norm E}
    \end{equation}
    \label{lem:bound of E_F}
\end{lemma}

By using the triangle inequality and 
norm sub-multiplicativity 
in \eqref{eq:a tilde}, and then 
substituting \eqref{eq:norm F} and \eqref{eq:norm E},
we obtain
\begin{align*}
  \|\mathbf{a^*} - \mathbf{\tilde{a}}\|_2 
  &\leq 
    \|E\|_2 \|\Phi^{\dagger}\|_2 \|\mathbf{a^*}\|_2 
    + \|E\|_2^2 \|\Phi^{\dagger}\|_2^2 \|\mathbf{a^*}\|_2 
    \nonumber \\
    &\qquad + \frac{
      \|E\|_2 D \|\Phi^{\dagger}\|_2^4 
    }{
      1 - D \|\Phi^{\dagger}\|_2^2
    }
    (\|\Phi\|_2 + \|E\|_2)
    \|\mathbf{a^*}\|_2
  \nonumber \\
  & \leq  
  \frac{
    4 e^{\frac{4}{\sigma^2}}\sqrt{M}
    \|\mathbf{a^*}\|_2
  }{
    \sigma^2 \sigma_{\min}(\Phi)
  } \|\mathbf{\tilde{t}}-\mathbf{t^*}\|_2
  + O(\|\mathbf{\tilde{t}}-\mathbf{t^*}\|_2^2),
\end{align*}
which is the bound given in Theorem \ref{thm:a dep t}.
Note that because 
$\|E\|_2 = O(\|\tilde{\mathbf{t}} -\mathbf{t^*}\|_2)$
(see \eqref{eq:norm E}), 
the first
term is the only term that is 
$O(\|\tilde{\mathbf{t}} -\mathbf{t^*}\|_2)$
in the first inequality above, so
the other terms are included in the 
$O(\|\tilde{\mathbf{t}} -\mathbf{t^*}\|_2^2)$ 
term at the end.

Lastly, we return to condition \eqref{eq:cond Neumann}, 
which must be satisfied in order for the 
bound above to hold. By using norm sub-multiplicativity
and the bound on $\|\Delta\|_2$ from \eqref{eq:def D}, 
we obtain
\begin{equation}
    \| 
        {\Phi^{\dagger}}^T \Delta \Phi^{\dagger} 
    \|_2
    \leq 
    \| \Phi^{\dagger}\|^2_2 \|E\|^2_2
    + 2 \|\Phi\|_2 \|\Phi^{\dagger}\|^2_2 \|E\|_2
    \label{eq:rhs cond}
\end{equation}
and by requiring that the right hand side above is
less than one, we obtain a quadratic constraint on $\|E\|_2$, 
satisfied if
\begin{equation*}
    \|E\|_2 <
    \sigma_{\max}(\Phi)\left(
        \sqrt{
            1 + \frac{
                \sigma^2_{\min}(\Phi)
            }{
                \sigma^2_{\max}(\Phi) 
            }
        } -1
    \right).
\end{equation*}
By using the bound on $\|E\|_2$ from \eqref{eq:norm E},
the above holds if \eqref{eq:cond t} holds.
Note that by imposing this, we also ensure that the condition
for the series in \eqref{eq:bound norm A} to converge holds, 
since $D\|\Phi^{\dagger}\|_2^2$ is equal to the right hand side
of \eqref{eq:rhs cond}.

\section{Conclusion}

In this paper, we proved primal stability in the non-negative
super-resolution problem, when addressed via convex duality. 
The main ingredient in our analysis is a quantitative version 
of the implicit function theorem, a folklore result in the 
theory of dynamical systems community. 

Our results provide precise orders in the number of measurements 
for the accuracy of the solution to the convex dual problem and 
an $\ell_\infty$ error bound on the primal spike locations.  

Future plans include the study of the dual approach to the noisy
super-resolution problem using similar techniques to the ones developed here.


\begin{thebibliography}{00}
	
\bibitem{candes2014towards}
    Emmanuel J. Cand{\`e}s and Carlos Fernandez-Granda,
    ``Towards a mathematical theory of super-resolution'',
    Communications on Pure and Applied Mathematics,
    67(6), 906--956 (2014)

\bibitem{schiebinger} 
    Geoffrey Schiebinger, Elina Robeva and Benjamin Recht,
    ``Superresolution without separation'',
    Information and Inference: A Journal of the IMA, 7(1), 1--30 (2018)

\bibitem{boyd2017}
    Nicholas Boyd, Geoffrey Schiebinger and Benjamin Recht,
    ``The alternating descent conditional gradient method for sparse inverse problems'',
    SIAM Journal on Optimization, 27(2), 616--639 (2017)

\bibitem{NesterovOpt} 
  Yurii Nesterov,
  ``Introductory Lectures on Convex Optimization: A Basic Course'', 
  Springer Publishing Company (2004)

\bibitem{Stewart90perturbationtheory} 
  G. W. Stewart,
  ``Perturbation theory and least squares with errors in the variables'', 
  Contemporary Mathematics 112: Statistical Analysis of Measurement 
  Error Models and Applications, 
  American Mathematical Society, 171--181 (1990)

\bibitem{lopez2007semiinfinite} 
  Marco L{\'o}pez and Georg Still, 
  ``Semi-infinite programming'', 
  European Journal of Operational Research, 180(2), 491--518 (2007)

\bibitem{eftekhari2019equivalence} 
  Armin Eftekhari and Andrew Thompson, 
  ``Sparse inverse problems over measures:
  equivalence of the conditional gradient and exchange methods'', 
  To appear in SIAM Journal on Optimization,
  arXiv: https://arxiv.org/abs/1804.10243
  (2019)

\bibitem{betzig2006imaging}
  Eric Betzig, George H. Patterson, Rachid Sougrat, O. Wolf Lindwasser, 
  Scott Olenych, Juan S. Bonifacino, Michael W. Davidson, 
  Jennifer Lippincott-Schwartz and Harald F. Hes
  ``Imaging intracellular fluorescent proteins at nanometer resolution'',
  Science, 313(5793), 1642--1645 (2006)

\bibitem{puschmann2005super}
  Klaus G. Puschmann and Franz Kneer,
  ``On super-resolution in astronomical imaging'',
  Astronomy \& Astrophysics, 436(1), 373--378 (2005)

\bibitem{tur2011innovation}
  Ronen Tur, Yonina C. Eldar and Zvi Friedman,
  ``Innovation rate sampling of pulse streams with application to 
  ultrasound imaging'',
  IEEE Transactions on Signal Processing, 59(4), 1827--1842 (2011)

\bibitem{liverani}
  Carlangelo Liverani,
  ``Implicit function theorem (a quantitative version)'',
  retrieved January 13, 2019, 
  from https://www.mat.uniroma2.it/~liverani/SysDyn15/app1.pdf

\end{thebibliography}
\end{document}